\documentclass[11pt,a4paper]{article}
\usepackage[T1]{fontenc}
\usepackage{mathpazo}
\usepackage[scaled=0.92]{helvet}
\usepackage{courier}
\usepackage[a4paper,top=2.6cm,bottom=2.8cm,left=2.9cm,right=2.9cm]{geometry}
\usepackage{amsmath,amssymb,amsthm}
\usepackage{microtype}
\usepackage{xcolor}
\usepackage{tikz}
\usetikzlibrary{arrows.meta,positioning,decorations.markings,bending}
\usepackage{url}

\usepackage[small,bf]{titlesec}
\titleformat{\section}{\normalfont\Large\bfseries}{\thesection}{0.8em}{}
\titleformat{\subsection}{\normalfont\large\bfseries}{\thesubsection}{0.8em}{}
\titlespacing*{\section}{0pt}{2.2ex plus 1ex minus .2ex}{1.1ex plus .2ex}
\titlespacing*{\subsection}{0pt}{1.8ex plus 1ex minus .2ex}{0.8ex plus .2ex}

\usepackage[font=small,labelfont=bf,labelsep=period]{caption}

\usepackage[hidelinks]{hyperref}

\makeatletter
\renewcommand{\maketitle}{%
  \begin{center}
    {\LARGE\bfseries Discrete Spoilers: Introducing Abstract\\[3pt] Undergraduate Mathematics Using Graph Theory\par}
    \vspace{1.2em}
    {\large Fernando Diaz-Diaz\par}
    \vspace{0.4em}
    {\small Department of Mathematics, Grupo Interdisciplinar de Sistemas Complejos (GISC)\\
    Universidad Carlos III de Madrid, 28911 Legan\'es, Spain\\
    \texttt{fddiaz@math.uc3m.es}\par}
  \end{center}
  \vspace{1em}}
\renewenvironment{abstract}{%
  \begin{list}{}{\leftmargin=2.2em\rightmargin=2.2em}\item\relax\small
  {\bfseries Abstract.}\ }{\end{list}}

\makeatother

\theoremstyle{plain}

\theoremstyle{definition}

\newcommand{\R}{\mathbb{R}}
\newcommand{\im}{\operatorname{im}}

\begin{document}

\maketitle

\begin{abstract}
Students often have serious difficulties with many key concepts of the undergraduate mathematics curriculum. We argue that many of these difficulties can be eased by explicitly separating the key ideas from the technical details. To do so, we present finite graphs as a setting where the basic concepts appear more transparently, because on a graph most proofs only require basic techniques from linear algebra. We illustrate this approach in three areas that cover a large part of the curriculum. First, we explain basic operators from vector calculus, such as the gradient and the divergence, through the incidence matrix of a graph. Second, we show that the diffusion equation has a simple solution on a graph. We also solve the graph Dirichlet problem, in which the values at some vertices are fixed. Third, we introduce a second incidence matrix that illustrates homology, cohomology and the Helmholtz decomposition in a simple way. Finally, we outline how to apply this approach to other branches of mathematics.
\end{abstract}

\vspace{4px}

\begin{list}{}{\leftmargin=2.2em\rightmargin=2.2em}\item\relax\small{\bfseries Keywords.} graph theory, Laplacian, discrete calculus, homology, mathematics education\end{list}

\vspace{4px}

\section{Introduction}

In 1750, in a letter to Goldbach, Euler stated that a polyhedron with $V$ vertices, $E$ edges and $F$ faces satisfies $V-E+F=2$ \cite{richeson2008}. This identity was ``proven'' by him soon after, and again by Cauchy in 1813. Yet that same year Lhuilier found polyhedra that did not satisfy the formula \cite{lhuilier1813}, such as a polyhedron with the shape of a picture frame, for which $V-E+F=0$ (Figure~\ref{fig:euler}). These counterexamples forced mathematicians to change the statement of the theorem, its proof and even the definition of a polyhedron. Lhuilier himself showed that a polyhedron with $g$ tunnels satisfies $V-E+F=2-2g$. Later, Poincar\'e generalized the alternating sum $V-E+F$ to the Euler characteristic of a general space and related it to the Betti numbers of the space \cite{richeson2008}. Overall, the full understanding of Euler's formula required more than a century of collective work.

This history is a good illustration of the way mathematics advances. Mathematicians start from concrete cases, such as the cube and other convex polyhedra, make a conjecture and look for a proof. Using counterexamples, they can also refine the hypotheses that the theorem needs. At the end of the process, once they understand the problem, they extract a definition that condenses what they have learned \cite{lakatos1976}. In many courses, however, the order is the opposite. For example, a course on algebraic topology may define the Euler characteristic of a space as ``the alternating sum of the dimensions of the quotients of the kernel of each boundary map by the image of the next one'' before giving any motivation or example. Do we really expect students to understand this chain of abstract concepts without any preparation? More likely, students either understand very little and become frustrated, or learn to manipulate abstract concepts without understanding their meaning or the problems they are meant to solve \cite{skemp1976}. In the latter case, they learn mathematics, but not \emph{how to think like a mathematician}.

Of course, a course cannot reproduce every attempt and every counterexample that led to a definition. However, we can show students \emph{why} we introduce a definition and \emph{which} problem it aims to solve, before we present its technical version. In this article we argue that finite graphs are a natural setting to explore many abstract concepts of the undergraduate curriculum in this way \cite{doyle1984,lim2020,knill2012,grady2010}. We see two reasons that make graphs an ideal didactic tool. First, graphs are accessible. We can draw a graph, and in a proof we can draw each step of the argument and see what it means on a particular graph. Hence a student can explore mathematical arguments with much less abstraction, in the same spirit as the exploration of concrete cases that we described above. Second, and most importantly, on a finite graph, the study of operators acting on functions reduces to the study of matrices acting on vectors. Hence the theorems that characterize these operators are based on linear algebra, a subject that students already know well, and the proofs are short and intuitive.

Therefore, if we first teach a problem in a discrete setting, the student can put all their mental effort into the problem itself. Why is the question interesting? How do we translate it into mathematical language? Which ideas lead to a proof? Only after the student has understood the object of study and the ideas behind it can we introduce the continuous version and the difficulties that require new mathematical tools (compactness, completeness, etc.). We can see a problem on a graph as a \emph{discrete spoiler} of the more general problem: the graph version gives away the key idea, the ``click'' in the student's mind, without all the development needed to reach it in the more difficult setting. Graph versions are also \emph{discreet} spoilers, because the student realizes only later, with the hard version of the problem, that they were foreshadowing something important.
\begin{figure}[t]
\centering
\begin{tikzpicture}[solid/.style={line width=0.7pt},hid/.style={line width=0.5pt,dash pattern=on 2pt off 1.5pt,black!45}]
  \draw[solid] (0.000,0.000) -- (0.000,1.500);
  \draw[hid] (0.001,0.002) -- (0.540,0.720);
  \draw[solid] (0.000,0.000) -- (1.500,0.000);
  \draw[solid] (0.000,1.500) -- (0.540,2.220);
  \draw[solid] (0.000,1.500) -- (1.500,1.500);
  \draw[hid] (0.540,0.720) -- (0.540,2.216);
  \draw[hid] (0.540,0.720) -- (2.036,0.720);
  \draw[solid] (0.540,2.220) -- (2.040,2.220);
  \draw[solid] (1.500,0.000) -- (1.500,1.500);
  \draw[solid] (1.500,0.000) -- (2.040,0.720);
  \draw[solid] (1.500,1.500) -- (2.040,2.220);
  \draw[solid] (2.040,0.720) -- (2.040,2.220);
  \draw[solid] (3.400,0.360) -- (6.000,0.360);
  \draw[solid] (3.400,0.000) -- (6.000,0.000);
  \draw[solid] (4.488,0.504) -- (5.488,0.504);
  \draw[hid] (4.488,0.384) -- (5.488,0.384);
  \draw[solid] (3.400,0.360) -- (4.488,0.504);
  \draw[hid] (3.400,0.000) -- (4.488,0.384);
  \draw[solid] (3.400,0.000) -- (3.400,0.360);
  \draw[hid] (4.488,0.384) -- (4.488,0.504);
  \draw[solid] (6.000,0.360) -- (6.936,1.608);
  \draw[solid] (6.000,0.000) -- (6.936,1.248);
  \draw[solid] (5.488,0.504) -- (5.848,0.984);
  \draw[hid] (5.488,0.384) -- (5.848,0.864);
  \draw[solid] (6.000,0.360) -- (5.488,0.504);
  \draw[hid] (6.000,0.000) -- (5.488,0.384);
  \draw[solid] (6.000,0.000) -- (6.000,0.360);
  \draw[hid] (5.488,0.384) -- (5.488,0.504);
  \draw[solid] (6.936,1.608) -- (4.336,1.608);
  \draw[hid] (6.936,1.248) -- (4.336,1.248);
  \draw[solid] (5.848,0.984) -- (4.848,0.984);
  \draw[hid] (5.848,0.864) -- (5.758,0.864);
  \draw[solid] (5.758,0.864) -- (4.848,0.864);
  \draw[solid] (6.936,1.608) -- (5.848,0.984);
  \draw[hid] (6.936,1.248) -- (5.848,0.864);
  \draw[solid] (6.936,1.248) -- (6.936,1.608);
  \draw[hid] (5.848,0.864) -- (5.848,0.984);
  \draw[solid] (4.336,1.608) -- (3.400,0.360);
  \draw[hid] (4.336,1.248) -- (3.400,0.000);
  \draw[solid] (4.848,0.984) -- (4.488,0.504);
  \draw[solid] (4.848,0.864) -- (4.578,0.504);
  \draw[hid] (4.578,0.504) -- (4.488,0.384);
  \draw[solid] (4.336,1.608) -- (4.848,0.984);
  \draw[hid] (4.336,1.248) -- (4.848,0.864);
  \draw[hid] (4.336,1.248) -- (4.336,1.608);
  \draw[solid] (4.848,0.864) -- (4.848,0.984);
  \node[font=\small] at (1.02,-0.45) {$V=8,\ E=12,\ F=6$};
  \node[font=\small] at (1.02,-0.95) {$8-12+6=2$};
  \node[font=\small] at (5.17,-0.45) {$V=16,\ E=32,\ F=16$};
  \node[font=\small] at (5.17,-0.95) {$16-32+16\neq2$};
\end{tikzpicture}
\caption{Left, a cube. Right, a polyhedron in the shape of a beveled picture frame, with a square tunnel through the middle. Below each polyhedron, the numbers $V$, $E$ and $F$ of vertices, edges and faces, and the value of $V-E+F$ compared with $2$. Dashed lines are hidden edges.}
\label{fig:euler}
\end{figure}
In the rest of this article, we exemplify this approach in three areas of the curriculum: vector calculus, partial differential equations and algebraic topology. In Section~\ref{sec:calculus}, we build the gradient and the divergence of a graph from its incidence matrix $B$, with a model of opinions in a social network as motivation, and study the four fundamental subspaces of $B$. In Section~\ref{sec:diffusion}, we solve this model, which is the graph version of the diffusion equation, and the analogous problem with zealots (people who never change their opinion), which is the graph version of the Dirichlet problem. In Section~\ref{sec:topology}, we build a curl by looking at the triangles of the graph and obtain homology, cohomology and the complete Helmholtz decomposition of the edge space. In Section~\ref{sec:other}, we outline how to apply the same approach to other branches of mathematics. We end in Section~\ref{sec:discussion} with a discussion of earlier work and of the limits of our approach.

\section{Calculus on a graph}\label{sec:calculus}

We start by motivating our framework with a question from everyday life. Several people discuss a topic, for example the quality of a film, with every person having an initial opinion on the film. Every day, each person talks only with his or her friends, and after each conversation the opinions of the two friends are closer than before. The pairs of friends form a graph $G=(V,E)$, with a set $V$ of vertices, one for each person, and a set $E$ of edges, one for each pair of friends. We will assume that $G$ is undirected, unweighted and simple (no loops or multiple edges). We write $n$ for the number of vertices and $m$ for the number of edges. How do the opinions change in time? Do all the people reach an agreement, and on which opinion?

To answer these questions, we need a mathematical description of the opinions and a rule for the way in which they change. We will assume that the opinion of the person at each vertex $v$ at time $t$ is a real number $x_v(t)$, with large numbers for a good opinion and small numbers for a bad one. Mathematically, this assumption corresponds, at each time $t$, to a function from $V$ to $\R$ that assigns the real number $x_v(t)$ to each vertex $v$. Since $V$ is finite, this function is determined by its $n$ values $x_v(t)$, so we can identify it with the vector $x(t)\in\R^n$ and the space of all opinions with $\R^n$. This finite description is very different from the case of functions from $\R$ to $\R$, which form a space of infinite dimension and require formal tools from real analysis, such as suprema or Cauchy sequences. Already at this point, we see that working on a graph greatly reduces the technical tools that we need to do mathematics.\footnote{In Section~\ref{sec:diffusion}, we will introduce time as a continuous variable, so technically we still need these tools to understand the derivatives with respect to $t$. However, the only equations in time that we meet are linear systems of ordinary differential equations with constant coefficients, and all the difficulties of the continuum that we discuss in this section come from spatial terms, not temporal.}

For the rule that governs the evolution, we will assume that the opinions change in the simplest possible way. People change their opinion only through conversations with their friends. In a conversation, each friend is influenced by the other in the same way, so the two opinions move toward each other by the same amount. In particular, the friend with the larger opinion lowers it, while the friend with the smaller opinion raises it. Two friends with very different opinions have more to argue about than two friends who almost agree, so their opinions change faster. To model these assumptions, we introduce a flow of opinion along each edge, from the larger opinion to the smaller one, proportional to the difference between the two opinions. In other words, the flow along each edge is proportional to a discrete version of the derivative of the opinions.

\begin{figure}[t]
\centering
\begin{tikzpicture}[line width=0.6pt,>={Stealth[length=2mm]},
  vtx/.style={circle,draw,fill=white,inner sep=1.5pt,minimum size=15pt,font=\small},
  lab/.style={font=\footnotesize}]
  \node[vtx] (1) at (0,0) {$1$};
  \node[vtx] (2) at (1.8,0) {$2$};
  \node[vtx] (3) at (0,1.8) {$3$};
  \node[vtx] (4) at (1.8,1.8) {$4$};
  \node[vtx] (5) at (0.9,2.9) {$5$};
  \draw[->] (1) -- node[lab,below] {$e_1$} (2);
  \draw[->] (1) -- node[lab,left] {$e_2$} (3);
  \draw[->] (2) -- node[lab,right] {$e_3$} (4);
  \draw[->] (3) -- node[lab,below] {$e_4$} (4);
  \draw[->] (3) -- node[lab,above left] {$e_5$} (5);
  \draw[->] (4) -- node[lab,above right] {$e_6$} (5);
\end{tikzpicture}
\caption{The house graph, with one person at each vertex and an edge between each pair of friends. Each edge is given an arbitrary orientation.}
\label{fig:house}
\end{figure}
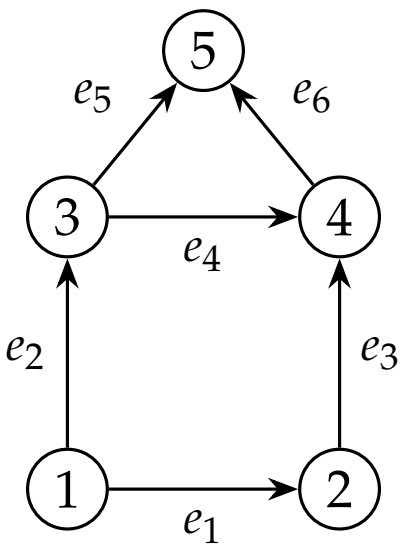

How can we define a discrete derivative? In one dimension, the derivative of a function $x$ at a point $v$ is the limit of $(x(v+\varepsilon)-x(v))/\varepsilon$ as $\varepsilon\to0$. But on a graph, this limit has no meaning, because two vertices cannot be infinitely close. The smallest distance between two vertices is one, so we set $\varepsilon=1$ and identify the point $v+\varepsilon$ with a neighbor of $v$, say, $w$. The derivative of $x$ along the edge $(v,w)$ is then the difference $x_w-x_v$. Note that the discrete derivative is a linear map between functions on two different sets, namely from functions on the vertices, in $\R^n$, to functions on the edges, in $\R^m$. In the continuum, the gradient also maps a function to a different kind of object, a vector field. However, the function and the vector field are defined on the same domain, so many students do not notice this change. On a graph, in contrast, this change is much more explicit, because the two functions are defined on different sets.

Our definition of the discrete derivative still has a problem. An undirected edge between $v$ and $w$ has no preferred endpoint, so we could equally take the difference $x_v-x_w$. The first difference corresponds to the orientation of the edge from $v$ to $w$, and the second one to the reverse. Which of the two should we choose? It turns out that the choice does not matter, as long as we use the same orientation throughout~\cite{biggs1993}, so we will choose an arbitrary orientation for each edge. Figure~\ref{fig:house} shows one such choice in an example with five people (we will refer to this graph as the house graph because of its shape). So, in summary, if we write $t(e)$ for the vertex where an edge $e$ starts, the tail of $e$, and $h(e)$ for the vertex where it ends, the head of $e$, then the discrete derivative of $x$ along $e$, which we call the gradient of $x$ along $e$ and denote by $(\nabla x)_e$, is
\[
(\nabla x)_e=x_{h(e)}-x_{t(e)}.
\]

The gradient is linear in $x$. On a graph, every linear operator on the functions on the vertices is a matrix, so we can write the gradient as the product of a matrix and the vector $x$. What does this matrix look like? First, it must clearly be rectangular, with $m$ rows and $n$ columns. Since each row $e$ computes a difference between  the two endpoints of $e$, it must have a value $1$ in column $h(e)$, $-1$ in column $t(e)$ and zeros elsewhere. For example, for the house graph with the orientation of Figure~\ref{fig:house}, the matrix of the gradient is
\begin{equation}\label{eq:grad}
\begin{pmatrix}
-1& 1& 0& 0& 0\\
-1& 0& 1& 0& 0\\
 0&-1& 0& 1& 0\\
 0& 0&-1& 1& 0\\
 0& 0&-1& 0& 1\\
 0& 0& 0&-1& 1
\end{pmatrix}.
\end{equation}

Now that we have a notion of gradient, we can define the flow of opinion as the vector $f=-\nabla x\in\R^m$. Each component $f_e$ is the amount of opinion that moves along the edge $e$ per unit of time, positive if opinion moves from $t(e)$ to $h(e)$ and negative otherwise. We choose this definition by analogy with Fourier's law, according to which the flow of heat is $F=-\nabla u$, with $u$ the temperature.

Our next goal is to combine the flows to deduce how the opinion at each vertex changes. The opinion at vertex $v$ changes at a rate equal to the net flow into $v$, which is
\[
\underbrace{\sum_{e:\,h(e)=v}f_e}_{\text{total flow into }v} - \underbrace{\sum_{e:\,t(e)=v}f_e}_{\text{total flow out of }v}.
\]
We identify this net flow with minus the divergence of $f$ at $v$, $-(\operatorname{div}f)_v$, by analogy with the continuum. Note that div maps functions on the edges to functions on the vertices. We can once again find the matrix associated with this linear operator. Each row $v$ has the entry $1$ in the columns of the edges that start at $v$, the entry $-1$ in the columns of the edges that end at $v$, and zeros elsewhere. For the house graph with the orientation of Figure~\ref{fig:house}, the matrix of the divergence is
\[
\begin{pmatrix}
 1& 1& 0& 0& 0& 0\\
-1& 0& 1& 0& 0& 0\\
 0&-1& 0& 1& 1& 0\\
 0& 0&-1&-1& 0& 1\\
 0& 0& 0& 0&-1&-1
\end{pmatrix}.
\]
But this matrix is minus the transpose of \eqref{eq:grad}! We defined the gradient as a discrete derivative and the divergence as a balance of flows, and yet the two operators have the same matrix up to a transpose and a sign. This identity is true for every graph, and as a result a single matrix is enough to describe both operators. We call the \emph{incidence matrix} of the graph the $n\times m$ matrix $B$ with $B_{ve}=1$ if $v=h(e)$, $B_{ve}=-1$ if $v=t(e)$ and $B_{ve}=0$ otherwise, so that $\nabla=B^T$ and $\operatorname{div}=-B$. We use the name \emph{incidence} because the nonzero entries of each column $e$ of $B$ indicate the two vertices incident to $e$.

\subsection{The four fundamental subspaces of $B$}\label{sec:subspaces}

Before we go back to our model of opinions, let us see what the matrix $B$ tells us about the graph. Like any other matrix, $B$ has four fundamental subspaces, namely the kernel and the image of $B$ and of its transpose \cite{strang1993}. Do these subspaces contain information about the graph? We will see that, indeed, they describe its connected components, its cycles and its cuts. 

We can find the kernel of $B^T$ directly. A vector $x\in\R^n$ satisfies $B^Tx=0$ if and only if it has the same value at the two endpoints of every edge. In that case, $x$ has the same value along every path, so it must be constant on each connected component. Conversely, every vector that is constant on each connected component satisfies $B^Tx=0$. Hence $\dim\ker B^T=c$, where $c$ is the number of connected components of the graph. In our model of opinions, the opinion vectors in $\ker B^T$ produce no flow. Once we know the kernel of $B^T$, we obtain the image of $B$ for free, as the two subspaces are orthogonal complements in $\R^n$, so $\im B$ consists of the vectors with zero sum on each connected component. This has a natural interpretation in our opinion model. $Bf$ represents the net change in opinions induced by a flow $f$. Each flow transfers opinion from one vertex to another without changing the component's total, so only opinion changes that preserve these totals are admissible. $\im B$ consists precisely of the admissible opinion changes.

Let us now look at the kernel of $B$. Since $B=-\operatorname{div}$, this kernel contains the flows with zero divergence, that is, the flows for which the total flow that enters each vertex is equal to the total flow that leaves it. The simplest way to obtain such a flow is to send one unit of flow around a cycle of the graph. We describe this flow with the vector $\gamma$$\in\R^m$ that is equal to $1$ on the edges of the cycle whose orientation agrees with the direction of the flow, to $-1$ on the other edges of the cycle and to $0$ elsewhere. By construction, this flow belongs to $\ker B$, and it can be proven that analogous flows around each cycle of the graph span the kernel of $B$~\cite{biggs1993}. From now on, we will call every vector $\gamma$$\in\R^m$ with $B\gamma=0$ a \emph{cycle}. For example, on the house graph, the flows around the square and around the roof are the cycles $\gamma_1$$=(1,-1,1,-1,0,0)^T$ and $\gamma_2$$=(0,0,0,1,-1,1)^T$. In our model, a cycle is a flow that does not change any opinion.

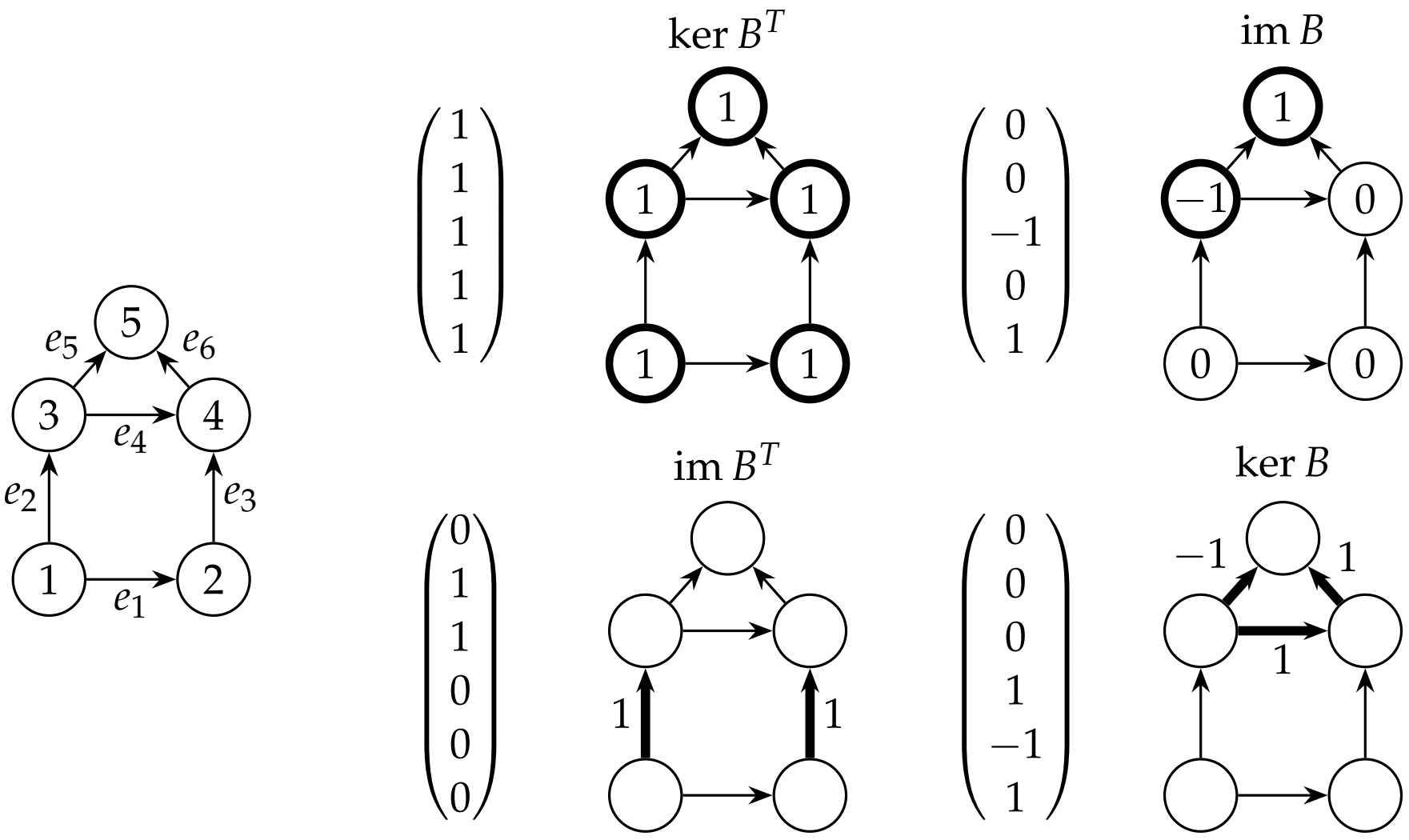
\begin{figure}[t]
\centering
\begin{tikzpicture}[x=0.85cm,y=0.85cm,line width=0.6pt,>={Stealth[length=2mm]},
  vtx/.style={circle,draw,inner sep=0.5pt,minimum size=17pt,font=\small},
  elab/.style={font=\small,inner sep=2pt}]
\begin{scope}
  \node[vtx] (n1) at (0.00,-2.10) {$1$};
  \node[vtx] (n2) at (1.60,-2.10) {$2$};
  \node[vtx] (n3) at (0.00,-0.50) {$3$};
  \node[vtx] (n4) at (1.60,-0.50) {$4$};
  \node[vtx] (n5) at (0.80,0.40) {$5$};
  \draw[->] (n1) -- node[elab,below] {$e_1$} (n2);
  \draw[->] (n1) -- node[elab,left] {$e_2$} (n3);
  \draw[->] (n2) -- node[elab,right] {$e_3$} (n4);
  \draw[->] (n3) -- node[elab,below] {$e_4$} (n4);
  \draw[->] (n3) -- node[elab,above left] {$e_5$} (n5);
  \draw[->] (n4) -- node[elab,above right] {$e_6$} (n5);
\end{scope}
\begin{scope}
  \node[vtx,line width=1.8pt] (n1) at (5.80,0.00) {$1$};
  \node[vtx,line width=1.8pt] (n2) at (7.40,0.00) {$1$};
  \node[vtx,line width=1.8pt] (n3) at (5.80,1.60) {$1$};
  \node[vtx,line width=1.8pt] (n4) at (7.40,1.60) {$1$};
  \node[vtx,line width=1.8pt] (n5) at (6.60,2.50) {$1$};
  \draw[->] (n1) -- (n2);
  \draw[->] (n1) -- (n3);
  \draw[->] (n2) -- (n4);
  \draw[->] (n3) -- (n4);
  \draw[->] (n3) -- (n5);
  \draw[->] (n4) -- (n5);
  \node[font=\small] at (6.60,3.25) {$\ker B^T$};
  \node[font=\small] at (4.0,1.25) {$\arraycolsep=5pt\begin{pmatrix}\,1\,\\\,1\,\\\,1\,\\\,1\,\\\,1\,\end{pmatrix}$};
\end{scope}
\begin{scope}
  \node[vtx] (n1) at (11.20,0.00) {$0$};
  \node[vtx] (n2) at (12.80,0.00) {$0$};
  \node[vtx,line width=1.8pt] (n3) at (11.20,1.60) {$-1$};
  \node[vtx] (n4) at (12.80,1.60) {$0$};
  \node[vtx,line width=1.8pt] (n5) at (12.00,2.50) {$1$};
  \draw[->] (n1) -- (n2);
  \draw[->] (n1) -- (n3);
  \draw[->] (n2) -- (n4);
  \draw[->] (n3) -- (n4);
  \draw[->] (n3) -- (n5);
  \draw[->] (n4) -- (n5);
  \node[font=\small] at (12.00,3.25) {$\operatorname{im}B$};
  \node[font=\small] at (9.4,1.25) {$\arraycolsep=5pt\begin{pmatrix}0\\0\\-1\\0\\1\end{pmatrix}$};
\end{scope}
\begin{scope}
  \node[vtx] (n1) at (5.80,-4.20) {};
  \node[vtx] (n2) at (7.40,-4.20) {};
  \node[vtx] (n3) at (5.80,-2.60) {};
  \node[vtx] (n4) at (7.40,-2.60) {};
  \node[vtx] (n5) at (6.60,-1.70) {};
  \draw[->] (n1) -- (n2);
  \draw[->,line width=2.2pt] (n1) -- node[elab,left] {$1$} (n3);
  \draw[->,line width=2.2pt] (n2) -- node[elab,right] {$1$} (n4);
  \draw[->] (n3) -- (n4);
  \draw[->] (n3) -- (n5);
  \draw[->] (n4) -- (n5);
  \node[font=\small] at (6.60,-0.95) {$\operatorname{im}B^T$};
  \node[font=\small] at (4.0,-2.95) {$\arraycolsep=5pt\begin{pmatrix}0\\1\\1\\0\\0\\0\end{pmatrix}$};
\end{scope}
\begin{scope}
  \node[vtx] (n1) at (11.20,-4.20) {};
  \node[vtx] (n2) at (12.80,-4.20) {};
  \node[vtx] (n3) at (11.20,-2.60) {};
  \node[vtx] (n4) at (12.80,-2.60) {};
  \node[vtx] (n5) at (12.00,-1.70) {};
  \draw[->] (n1) -- (n2);
  \draw[->] (n1) -- (n3);
  \draw[->] (n2) -- (n4);
  \draw[->,line width=2.2pt] (n3) -- node[elab,below] {$1$} (n4);
  \draw[->,line width=2.2pt] (n3) -- node[elab,above left] {$-1$} (n5);
  \draw[->,line width=2.2pt] (n4) -- node[elab,above right] {$1$} (n5);
  \node[font=\small] at (12.00,-0.95) {$\ker B$};
  \node[font=\small] at (9.4,-2.95) {$\arraycolsep=5pt\begin{pmatrix}0\\0\\0\\1\\-1\\1\end{pmatrix}$};
\end{scope}
\end{tikzpicture}
\caption{Illustration of the four fundamental subspaces of the incidence matrix $B$ of the house graph, with the numbering and orientation of the left panel. For each subspace, we show one vector, both as a column vector and on the graph. For vectors of the vertex space $\R^n$, we write each entry inside its vertex and draw a thick border around the vertices with nonzero entries. For vectors of the edge space $\R^m$, we draw the edges with nonzero entries with thick lines and write these entries next to them. Top left, a constant opinion vector, in $\ker B^T$. Top right, the change of the opinions produced by one unit of flow along $e_5$, in $\operatorname{im}B$. Bottom left, the cut that separates the vertices $1$ and $2$ from the other vertices, in $\operatorname{im}B^T$, which is also the flow produced by the opinion vector equal to $1$ on the vertices $1$ and $2$ and to $0$ on the other vertices. Bottom right, a flow with zero divergence, that is, a cycle, in $\ker B$. In the two bottom panels, the flow on an edge is positive if it goes in the direction of the orientation of the edge in the left panel, and negative otherwise.}
\label{fig:subspaces}
\end{figure}

Finally, let us look at the image of $B^T$. To construct vectors in $\im B^T$, we can choose opinion vectors $x$ and apply $B^T$ to them. Consider the indicator vector $x$ of a subset $A$, which takes the value $1$ on the vertices in $A$ and $0$ on the vertices in its complement $A^c$.  With this choice, an edge can have nonzero flow only if it connects $A$ to $A^c$. The set of edges connecting a vertex subset to its complement is called a \emph{cut} in graph theory. Since every cut can be obtained from such a choice of $x$, the space $\im B^T$ is spanned by all cut vectors, that is, by the vectors $B^Tx$ with $x$ an indicator vector. Figure~\ref{fig:subspaces} shows the cut of $A=\{1,2\}$ in the house graph, together with one vector in each of the other three subspaces.

We can now put all these pieces together to see how any flow on a graph can be decomposed. Since $\operatorname{im}B^T$ and $\ker B$ are orthogonal complements in $\mathbb{R}^m$, every flow is, in a unique way, the sum of a gradient and a cycle. This is the graph version of the Helmholtz decomposition of a vector field into a gradient and a divergence-free field \cite{lim2020}. Students of vector calculus often meet the Helmholtz decomposition without a proof, since it requires tools that are typically introduced later in the curriculum. Here, in contrast, we obtain it directly from a basic result in linear algebra, namely that every vector decomposes uniquely into components belonging to two orthogonal complementary subspaces. 

\section{Diffusion dynamics on a graph}\label{sec:diffusion}

Let us now go back to our model of opinions. In this section, we will assume that $G$ is connected. In Section~\ref{sec:calculus}, we saw that the rate of change of the opinions is minus the divergence of the flow, which we can now write as $dx/dt=-\operatorname{div}f=Bf$. If we substitute the definition of the flow, $f=-\nabla x=-B^Tx$, we obtain
\begin{equation}\label{eq:dynamics}
\frac{dx}{dt}=-BB^Tx.
\end{equation}
We call the $n\times n$ matrix $L=BB^T$ the \emph{Laplacian} of the graph, so that \eqref{eq:dynamics} becomes the diffusion equation $dx/dt=-Lx$. We choose these names to mimic the continuum, where the operator $\operatorname{div}\circ\nabla$ is the Laplace operator $\Delta$ and the equation $\partial u/\partial t=\Delta u$ is the diffusion equation. $L$ is clearly symmetric, because $L^T=(BB^T)^T=BB^T=L$, and positive semidefinite, because $x^TLx=(B^Tx)^T(B^Tx)=\|B^Tx\|^2\ge0$ for every $x\in\R^n$. We call the quantity
\[
\mathcal{E}(x)=x^TLx=\|B^Tx\|^2=\sum_{e}\bigl(x_{h(e)}-x_{t(e)}\bigr)^2
\]
the energy of $x$, by analogy with the elastic energy of a network of springs placed on the edges. In our model, $\mathcal{E}(x)$ measures the total disagreement of the opinions $x$. It vanishes if and only if $B^Tx=0$, that is, if and only if $x\in\ker B^T$. Since $G$ is connected, $\ker B^T$ is spanned by the vector $\mathbf 1\in\R^n$ with all entries equal to one. This implies that $Lx=0$ only for the multiples of $\mathbf 1$, so $0$ is an eigenvalue of $L$ and its only eigenvectors are the multiples of $\mathbf 1$. Moreover, the total opinion $\sum_vx_v=\mathbf 1^Tx$ does not change in time, because $dx/dt=Bf$ belongs to the image of $B$, whose vectors have zero sum.

The solution of the linear system of ODEs \eqref{eq:dynamics} is $x(t)=e^{-tL}x(0)$, with $e^{-tL}=\sum_{k=0}^\infty(-tL)^k/k!$. To analyze this exponential, we will write $L$ in an orthonormal basis of eigenvectors $\varphi_1,\dots,\varphi_n$, with eigenvalues $\lambda_1\le\dots\le\lambda_n$. Because $L$ is positive semidefinite, all these eigenvalues are nonnegative. We expand the vector of initial opinions in this basis, $x(0)=\sum_kc_k\varphi_k$, with $c_k=\langle x(0),\varphi_k\rangle$ and $\langle\cdot,\cdot\rangle$ the usual inner product. Since $e^{-tL}\varphi_k=e^{-\lambda_kt}\varphi_k$, we obtain
\[
x(t)=e^{-tL}x(0)=\sum_{k=1}^nc_k\,e^{-\lambda_kt}\varphi_k .
\]
Each component with $\lambda_k>0$ decays exponentially at a rate equal to its eigenvalue. Since $G$ is connected, the only zero eigenvalue is $\lambda_1$, with $\varphi_1=\mathbf 1/\sqrt n$. Hence, for large $t$, $x(t)$ converges to $c_1\varphi_1$. In plain words, the opinion dynamics reach a steady state in which everybody agrees. This answers the question of Section~\ref{sec:calculus} about whether all the people reach an agreement. As for the opinion on which they agree, every opinion converges to $c_1/\sqrt n=\frac1n\sum_vx_v(0)$, the average of the initial opinions. 

Before we continue, let us look at these results from the point of view of a course on partial differential equations. While solving the model, we have derived with very little effort the graph versions of four important properties: the gradient and minus the divergence are adjoints of one another, the Laplacian is a positive semidefinite operator, the total opinion is conserved and the opinions converge to their average. All these facts follow from the factorization $L=BB^T$ and from basic results of linear algebra. In contrast, the analog results in the continuum require tools and derivations that are far from elementary. For example, to prove that the gradient and minus the divergence are adjoint, that is, that $\int_\Omega\nabla u\cdot F=-\int_\Omega u\operatorname{div}F$ for a function $u$ and a vector field $F$ on a domain $\Omega\subset\R^d$, we need to assume that $\Omega$ is bounded with a smooth boundary $\partial\Omega$, that $u$ and $F$ are smooth enough, and that $u$ vanishes on $\partial\Omega$. Similarly, proving that $-\Delta$, with zero boundary values, has a basis of eigenfunctions with positive eigenvalues requires tools from functional analysis, such as the characterization of compact self-adjoint operators \cite{evans2010}. We hope the reader is thus convinced that the difficulties in the continuum come from the technicalities of analysis, not from the diffusion process itself.

\subsection{Zealots and the Dirichlet problem}

Let us now change the model. Suppose that some people, the \textit{zealots}, never change their opinion. A zealot still talks with his or her friends and influences them. However, nobody influences the opinion of a zealot. We write $S$ for the set of vertices of the zealots and $I$ for the set of the other, influenceable, vertices. The opinion at each vertex $s\in S$ is a fixed number $z_s$, so the model becomes
\[
\frac{dx_v}{dt}=-(Lx)_v\quad\text{for } v\in I,\qquad x_s=z_s\quad\text{for } s\in S .
\]
We will focus on the stationary states of this model, which fulfill $dx_v/dt=0$ for every $v\in I$, that is,
\begin{equation}\label{eq:dirichlet}
x_s=z_s\quad\text{for } s\in S,\qquad (Lx)_v=0\quad\text{for } v\in I .
\end{equation}
In the continuum, \eqref{eq:dirichlet} is the Dirichlet problem, in which we look for a function $u$ with $\Delta u=0$ in a domain $\Omega$ and given values on its boundary $\partial\Omega$. The zealots are the graph version of the boundary \cite{doyle1984}.

To interpret the condition $(Lx)_v=0$ in \eqref{eq:dirichlet}, let us see how the operator $L$ acts on the opinion vector $x$. Using the definitions of discrete gradient and divergence, we obtain
\[
(Lx)_v=\sum_{e:\,h(e)=v}\bigl(x_v-x_{t(e)}\bigr)-\sum_{e:\,t(e)=v}\bigl(x_{h(e)}-x_v\bigr)=\sum_{w\sim v}(x_v-x_w),
\]
with $w\sim v$ meaning that $v$ and $w$ are neighbors. Denoting by $d_v$ the degree of $v$ and $\bar x_v$ the average opinion of its neighbors, we can express $(Lx)_v$ as
\begin{equation}\label{eq:mean}
(Lx)_v=d_v\,\bigl(x_v-\bar x_v\bigr).
\end{equation}

Hence, the condition $(Lx)_v=0$ states that the opinion at each vertex is the average of the opinions at the neighbors of $v$. We say that $x$ is \emph{harmonic} at $v$ if $(Lx)_v=0$.

Ordering the vertices so that those of $I$ come first, we can write the linear system \eqref{eq:dirichlet} in matrix form as
\begin{align} \label{eq:system_zealots}
\begin{pmatrix} L_{II} & L_{IS}\\ 0 & \mathrm{Id}\end{pmatrix}\begin{pmatrix} x_I\\ x_S\end{pmatrix}=\begin{pmatrix} 0\\ z\end{pmatrix},
\end{align}
with $L_{II}$ and $L_{IS}$ the submatrices of $L$ with rows in $I$ and columns in $I$ and in $S$, respectively, $x_I$ and $x_S$ the vectors of the opinions in $I$ and in $S$, $z$ the vector of the opinions of the zealots, and $\mathrm{Id}$ the identity matrix. The unknowns $x_I$ satisfy the system $L_{II}\,x_I=-L_{IS}\,z$, in which the opinions of the zealots act as the inhomogeneous term. The structure of the solutions depends only on the invertibility of $L_{II}$. Without zealots, $L_{II}=L$ is singular and the system has infinitely many solutions, corresponding to the multiples of $\mathbf 1$ discussed earlier in this section. With zealots, in contrast, $L_{II}$ is invertible and the system has exactly one solution. To prove this, let us show that the kernel of $L_{II}$ is trivial. Let $y\in\R^{|I|}$ satisfy $L_{II}\,y=0$, and let $\tilde y\in\R^n$ be equal to $y$ on $I$ and to $0$ on $S$. Then $\mathcal{E}(\tilde y)$$=\tilde y^TL\tilde y=y^TL_{II}\,y=0$, so $\tilde y$ is constant, and since $\tilde y$ vanishes on the nonempty set $S$, we obtain $\tilde y=0$ and therefore $y=0$.\footnote{This argument uses that $G$ is connected and that $S$ is nonempty. If $G$ were not connected, the people in a connected component without zealots could agree on any opinion, so the stationary state would not be unique.} Since $L_{II}$ is square, it is then invertible, and the system \eqref{eq:system_zealots} has exactly one solution, $x_I=-L_{II}^{-1}L_{IS} z$.

We have shown with relative ease that the Dirichlet problem on a graph has exactly one solution. In the continuum, uniqueness is still simple to prove, but existence is notoriously difficult. The reason is that, on a graph, we implicitly used uniqueness (the kernel of $L_{II}$ is trivial) to obtain existence (the system has a solution for any inhomogeneous term, i.e., $L_{II}$ is a surjective map). But in infinite dimension, a linear map with trivial kernel need not be surjective.\footnote{One example of this is the shift $(a_1,a_2,\dots)\mapsto(0,a_1,a_2,\dots)$ on infinite sequences of real numbers.} A course on partial differential equations therefore must prove existence with other arguments. The standard method obtains the solution as the minimizer of the energy $\int_\Omega|\nabla u|^2$ in a Hilbert space of functions with finite energy, and then uses regularity theory to show that this minimizer is a smooth harmonic function \cite{evans2010}.

We end this section with a historical note. The idea of minimizing the energy goes back to Riemann, who in 1851 used it to prove the existence of a solution of the Dirichlet problem. In the proof, he assumed that some function with the given boundary values attains the minimum energy. But in 1870, Weierstrass gave an example of an energy with an infimum that no function attains, showing that Riemann's proof had a gap. This gap was closed by Hilbert in 1900 \cite{monna1975}. Hence the existence of a minimizer, which on a graph can be proven by a first-year undergraduate, occupied Riemann, Weierstrass and Hilbert for half a century. Perhaps Riemann should have started with a graph.

\section{Topology on a graph}\label{sec:topology}

In Section~\ref{sec:calculus}, we saw how the four fundamental subspaces of the incidence matrix $B$ describe the connected components, the cycles and the cuts of a graph. Less obviously, we can also use incidence matrices to describe the holes of a graph. The study of holes, called homology, is a central topic of algebraic topology. In the continuum, homology requires difficult tools like chain complexes \cite{hatcher2002}. However, if we formulate the problem on a graph, we can once again derive the key results using only linear algebra. The price to pay is that we need a new type of incidence matrix.

Our new matrix can be derived from the third basic operator of vector calculus, the curl. In the continuum, the curl measures how much a field circulates around infinitesimal closed paths. On a graph we cannot build infinitesimal loops, so we use the closest proxies for them: triangles. The curl of a flow should therefore return its circulation around each triangle of the graph.\footnote{\label{fn:curl}This is another case in which the continuum hides the true nature of an operator. The curl at a point measures the circulation of the field around infinitesimally small closed paths. Since each of these paths lies in a plane, in $\R^3$ we can describe its circulation with the normal vector to the plane. In $\R^4$, however, a plane does not have a unique normal vector, so we cannot describe circulations with vectors. Hence the curl in $\R^4$ is not a vector field. Because vector calculus courses work almost always in $\R^3$, many students incorrectly assume that we can identify circulations with vectors in any dimension. The modern way to study the curl in any dimension is to understand it as a differential $2$-form.} To compute this circulation, we give each triangle an arbitrary orientation (clockwise or counterclockwise), as we did with the edges, and add the flows on its three edges, with a plus sign if the orientation of the edge agrees with that of the triangle and a minus sign otherwise. As with the gradient and the divergence, we can write this operator as a matrix, which we will denote by $B_2^T\in \R^{k \times m}$, where $k$ is the number of triangles in the graph. Each row of this matrix indicates the three edges of one triangle, with the entry $1$ or $-1$ depending on whether the orientation of the edge agrees with that of the triangle, and $0$ elsewhere. For example, the house graph has a single triangle, the roof, with the edges $e_4$, $e_5$ and $e_6$. If we orient it as $3\to4\to5\to3$, then
\[
B_2^T=\begin{pmatrix}0&0&0&1&-1&1\end{pmatrix}.
\]
On the other hand, the transpose $B_2$ maps a vector $w\in\R^k$, with one circulation per triangle, to the flow $B_2w$ that these circulations produce on the edges. To distinguish $B_2$ from the incidence matrix $B$ of Section~\ref{sec:calculus}, from now on we will write $B_1$ for the latter.

What happens if we multiply $B_1$ and $B_2$? To answer this question, let us take a vector of circulations $w\in\R^k$ and the flow $f=B_2w$ that it produces on the edges. $B_1f$ then captures minus the divergence, i.e. the net amount of the flow $f$ that accumulates at the vertex. But a flow that only circulates around triangles cannot accumulate at any vertex (all the flow that enters a vertex through one of its edges must leave it through the next edge of the triangle). Hence $B_1B_2w=B_1f=0$ for every $w\in\R^k$, so
\begin{equation}\label{eq:bb}
B_1B_2=0 .
\end{equation}
There is a profound interpretation of this identity. We can think of the two endpoints of an edge as its boundary. Similarly, the three edges of a triangle form its boundary. In these terms, the identity $B_1B_2=0$ says that the boundary of the boundary of a triangle is zero. This statement is the graph version of one of the central identities of advanced mathematics, ``the boundary of a boundary is zero'', which is the starting point of homology.

The identity $B_1B_2=0$ also connects with seemingly unrelated identities from vector calculus. First, as we saw above, $B_1B_2=0$ means that a circulation cannot accumulate flow at the vertices. Hence it is the graph version of the identity $\operatorname{div}\circ\operatorname{curl}=0$. Second, if we take transposes, we obtain $B_2^TB_1^T=0$. Since $B_1^T$ is the gradient and $B_2^T$ is the curl, this transposed identity is the graph version of $\operatorname{curl}\circ\nabla=0$. In terms of the fundamental subspaces, the two identities are the inclusions $\operatorname{im}B_2\subseteq\ker B_1$ (every flow induced by circulations has zero divergence) and $\operatorname{im}B_1^T\subseteq\ker B_2^T$ (every gradient flow has zero curl).

In vector calculus on $\R^3$, it is well known that the converse statements are also true: every field with zero curl is a gradient, and every field with zero divergence is the curl of another field. On a graph, however, the reverse inclusions can fail. For example, on the house graph, the cycle $\gamma_1$ around the square has zero divergence but is not in $\operatorname{im}B_2$ (every flow of the form $B_2w$ is zero on the edges that do not belong to a triangle, while $\gamma_1$ is nonzero on $e_1$, $e_2$ and $e_3$). This example shows that $\operatorname{im}B_2$ can be smaller than $\ker B_1$, so there are nonzero cycles in $\ker B_1\cap(\operatorname{im}B_2)^\perp$. Using the property $(\operatorname{im}B_2)^\perp=\ker B_2^T$, we can express this space as

\begin{equation}\label{eq:H}
H=\ker B_1\cap\ker B_2^T .
\end{equation}
Equation~\eqref{eq:H} highlights that the flows in $H$ have simultaneously zero divergence and zero curl. We will call these flows \emph{harmonic}. At first glance, the existence of these flows looks paradoxical. How can a flow neither accumulate at a vertex (zero divergence) nor go in circles (zero curl)? The way out of this paradox is that the curl $B_2^T$ only measures the circulation around triangles. A harmonic flow has zero circulation there, but it may circulate around longer cycles. What's more, we can exploit these flows to extract structural information about the graph, like the location of its holes. For example, on the house graph, the harmonic flows (proportional to $\gamma_1+\gamma_2/3$) go around the square, the only hole, but split at the roof to force the circulation around the triangle to vanish.

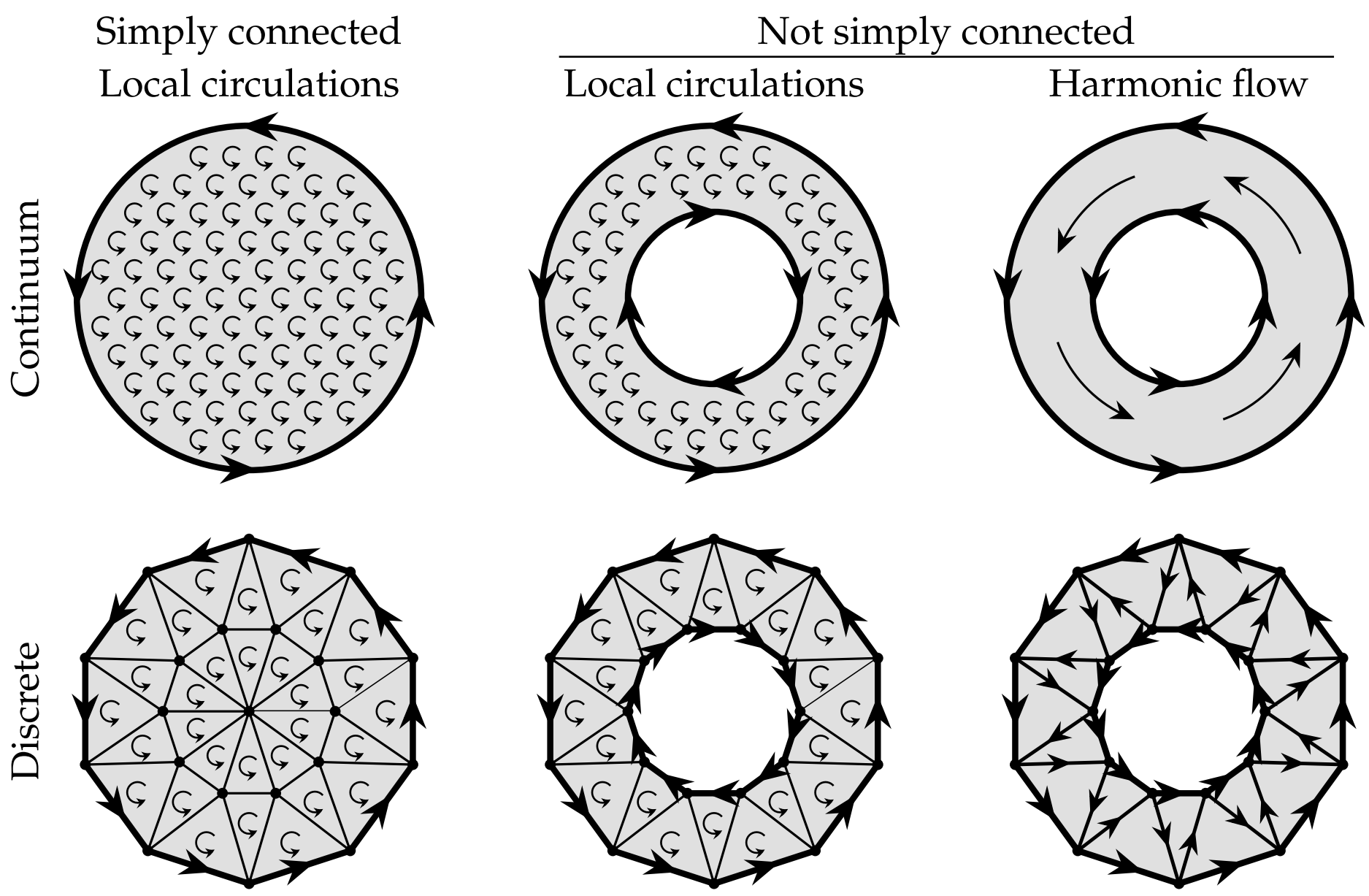
\begin{figure}[t]
\centering
\begin{tikzpicture}[x=0.78cm,y=0.78cm,line width=0.6pt,>={Stealth[length=1.6mm]},
  mid/.style={postaction={decorate},decoration={markings,mark=at position 0.65 with {\arrow{Stealth[length=3mm,width=2.4mm]}}}},
  midh/.style={postaction={decorate},decoration={markings,mark=at position 0.6 with {\arrow{Stealth[length=2.2mm,width=1.8mm]}}}},
  sl/.style={-{Stealth[length=0.9mm,width=0.9mm,bend]},line width=0.5pt}]
\begin{scope}
  \fill[black!12] (0,0) circle (2);
  \foreach \x/\y in {-1.71/-0.329,-1.71/0.329,-1.33/-0.987,-1.52/-0.658,-1.33/-0.329,-1.52/0.0,-1.33/0.329,-1.52/0.658,-1.33/0.987,-1.14/-1.316,-0.95/-0.987,-1.14/-0.658,-0.95/-0.329,-1.14/0.0,-0.95/0.329,-1.14/0.658,-0.95/0.987,-1.14/1.316,-0.57/-1.645,-0.76/-1.316,-0.57/-0.987,-0.76/-0.658,-0.57/-0.329,-0.76/0.0,-0.57/0.329,-0.76/0.658,-0.57/0.987,-0.76/1.316,-0.57/1.645,-0.19/-1.645,-0.38/-1.316,-0.19/-0.987,-0.38/-0.658,-0.19/-0.329,-0.38/0.0,-0.19/0.329,-0.38/0.658,-0.19/0.987,-0.38/1.316,-0.19/1.645,0.19/-1.645,0.0/-1.316,0.19/-0.987,0.0/-0.658,0.19/-0.329,0.0/0.0,0.19/0.329,0.0/0.658,0.19/0.987,0.0/1.316,0.19/1.645,0.57/-1.645,0.38/-1.316,0.57/-0.987,0.38/-0.658,0.57/-0.329,0.38/0.0,0.57/0.329,0.38/0.658,0.57/0.987,0.38/1.316,0.57/1.645,0.76/-1.316,0.95/-0.987,0.76/-0.658,0.95/-0.329,0.76/0.0,0.95/0.329,0.76/0.658,0.95/0.987,0.76/1.316,1.14/-1.316,1.33/-0.987,1.14/-0.658,1.33/-0.329,1.14/0.0,1.33/0.329,1.14/0.658,1.33/0.987,1.14/1.316,1.52/-0.658,1.71/-0.329,1.52/0.0,1.71/0.329,1.52/0.658}{\draw[sl] ([shift={(40:0.11)}]\x,\y) arc (40:320:0.11);}
  \draw[line width=1.6pt] (0,0) circle (2);
  \foreach \t in {0,90,180,270}{\draw[line width=1.6pt,-{Stealth[length=3.4mm,width=2.6mm]}] (\t-3:2) arc (\t-3:\t+3:2);}
\end{scope}
\begin{scope}[shift={(5.4,0)}]
  \fill[black!12,even odd rule] (0,0) circle (2) (0,0) circle (1);
  \foreach \x/\y in {-1.71/-0.329,-1.71/0.329,-1.33/-0.987,-1.52/-0.658,-1.33/-0.329,-1.52/0.0,-1.33/0.329,-1.52/0.658,-1.33/0.987,-1.14/-1.316,-0.95/-0.987,-1.14/-0.658,-1.14/0.658,-0.95/0.987,-1.14/1.316,-0.57/-1.645,-0.76/-1.316,-0.76/1.316,-0.57/1.645,-0.19/-1.645,-0.38/-1.316,-0.38/1.316,-0.19/1.645,0.19/-1.645,0.0/-1.316,0.0/1.316,0.19/1.645,0.57/-1.645,0.38/-1.316,0.38/1.316,0.57/1.645,0.76/-1.316,0.95/-0.987,0.95/0.987,0.76/1.316,1.14/-1.316,1.33/-0.987,1.14/-0.658,1.33/-0.329,1.33/0.329,1.14/0.658,1.33/0.987,1.14/1.316,1.52/-0.658,1.71/-0.329,1.52/0.0,1.71/0.329,1.52/0.658}{\draw[sl] ([shift={(40:0.11)}]\x,\y) arc (40:320:0.11);}
  \draw[line width=1.6pt] (0,0) circle (2);
  \draw[line width=1.6pt] (0,0) circle (1);
  \foreach \t in {0,90,180,270}{
    \draw[line width=1.6pt,-{Stealth[length=3.4mm,width=2.6mm]}] (\t-3:2) arc (\t-3:\t+3:2);
    \draw[line width=1.6pt,-{Stealth[length=3.4mm,width=2.6mm]}] (\t+5:1) arc (\t+5:\t-5:1);
  }
\end{scope}
\begin{scope}[shift={(10.8,0)}]
  \fill[black!12,even odd rule] (0,0) circle (2) (0,0) circle (1);
  \draw[line width=1.6pt] (0,0) circle (2);
  \draw[line width=1.6pt] (0,0) circle (1);
  \foreach \t in {0,90,180,270}{
    \draw[line width=1.6pt,-{Stealth[length=3.4mm,width=2.6mm]}] (\t-3:2) arc (\t-3:\t+3:2);
    \draw[line width=1.6pt,-{Stealth[length=3.4mm,width=2.6mm]}] (\t-5:1) arc (\t-5:\t+5:1);
  }
  \foreach \t in {45,135,225,315}{
    \draw[line width=0.6pt,-{Stealth[length=2.2mm,width=1.8mm]}] (\t-25:1.5) arc (\t-25:\t+25:1.5);
  }
\end{scope}
\begin{scope}[shift={(10.8,-4.8)}]
  \foreach \i in {0,...,9}{
    \pgfmathsetmacro{\a}{36*\i}\pgfmathsetmacro{\b}{36*(\i+1)}\pgfmathsetmacro{\c}{36*\i+18}\pgfmathsetmacro{\d}{36*(\i+1)+18}
    \fill[black!12] (\a:1) -- (\b:1) -- (\c:2) -- cycle;
    \fill[black!12] (\c:2) -- (\d:2) -- (\b:1) -- cycle;
  }
  \foreach \i in {0,...,9}{
    \pgfmathsetmacro{\a}{36*\i}\pgfmathsetmacro{\b}{36*(\i+1)}\pgfmathsetmacro{\c}{36*\i+18}\pgfmathsetmacro{\d}{36*(\i+1)+18}
    \draw[line width=1.6pt,mid] (\c:2) -- (\d:2);
    \draw[line width=1.6pt,mid] (\a:1) -- (\b:1);
    \draw[line width=0.9pt,midh] (\a:1) -- (\c:2);
    \draw[line width=0.9pt,midh] (\c:2) -- (\b:1);
    \fill (\a:1) circle (1.4pt); \fill (\c:2) circle (1.4pt);
  }
\end{scope}
\node[font=\small] at (0,3.05) {Simply connected};
\node[font=\small] at (8.1,3.05) {Not simply connected};
\draw[line width=0.5pt] (3.6,2.8) -- (12.6,2.8);
\node[font=\small] at (0,2.5) {Local circulations};
\node[font=\small] at (5.4,2.5) {Local circulations};
\node[font=\small] at (10.8,2.5) {Harmonic flow};
\node[font=\small,rotate=90] at (-2.6,0) {Continuum};
\node[font=\small,rotate=90] at (-2.6,-4.8) {Discrete};
\foreach \sx/\disk in {0/1,5.4/0}{
\begin{scope}[shift={(\sx,-4.8)}]
  \foreach \i in {0,...,9}{
    \pgfmathsetmacro{\a}{36*\i}\pgfmathsetmacro{\b}{36*(\i+1)}\pgfmathsetmacro{\c}{36*\i+18}\pgfmathsetmacro{\d}{36*(\i+1)+18}
    \fill[black!12] (\a:1) -- (\b:1) -- (\c:2) -- cycle;
    \fill[black!12] (\c:2) -- (\d:2) -- (\b:1) -- cycle;
    \ifnum\disk=1 \fill[black!12] (0,0) -- (\a:1) -- (\b:1) -- cycle; \draw (0,0) -- (\a:1); \fi
    \draw (\a:1) -- (\c:2) -- (\b:1);
    \pgfmathsetmacro{\ux}{(cos(\a)+cos(\b)+2*cos(\c))/3}\pgfmathsetmacro{\uy}{(sin(\a)+sin(\b)+2*sin(\c))/3}
    \pgfmathsetmacro{\vx}{(2*cos(\c)+2*cos(\d)+cos(\b))/3}\pgfmathsetmacro{\vy}{(2*sin(\c)+2*sin(\d)+sin(\b))/3}
    \draw[sl] ([shift={(40:0.12)}]\ux,\uy) arc (40:320:0.12);
    \draw[sl] ([shift={(40:0.12)}]\vx,\vy) arc (40:320:0.12);
    \ifnum\disk=1
      \pgfmathsetmacro{\wx}{(cos(\a)+cos(\b))/3}\pgfmathsetmacro{\wy}{(sin(\a)+sin(\b))/3}
      \draw[sl] ([shift={(40:0.11)}]\wx,\wy) arc (40:320:0.11);
      \draw (\a:1) -- (\b:1);
    \fi
  }
  \foreach \i in {0,...,9}{
    \pgfmathsetmacro{\a}{36*\i}\pgfmathsetmacro{\b}{36*(\i+1)}\pgfmathsetmacro{\c}{36*\i+18}\pgfmathsetmacro{\d}{36*(\i+1)+18}
    \draw[line width=1.6pt,mid] (\c:2) -- (\d:2);
    \ifnum\disk=0 \draw[line width=1.6pt,mid] (\b:1) -- (\a:1); \fi
    \fill (\a:1) circle (1.4pt); \fill (\c:2) circle (1.4pt);
  }
  \ifnum\disk=1 \fill (0,0) circle (1.4pt); \fi
\end{scope}
}
\end{tikzpicture}
\caption{Top, a disk (left) and a ring (middle) in the plane, with local counterclockwise circulations inside them, and the same ring with a flow that goes around the hole without local circulations (right). Bottom, triangulations of the disk and of the ring, with a counterclockwise circulation around each triangle (left and middle), and the harmonic flow of the triangulated ring (right), with thinner lines on the edges where its value is half. In the left and middle panels, the thick cycles are the global circulations that result from adding the local circulations.}
\label{fig:holes}
\end{figure}

A skeptical reader might doubt that this notion of hole captures the same idea as the empty center of a ring. However, mathematicians as early as Riemann have used the continuum analog of this definition to characterize holes \cite{riemann1857}: a surface has a hole if we can find a closed curve on it that does not bound any part of the surface, that is, a closed curve that we cannot write as a sum of small loops within the region enclosed by the curve. To better understand why such a curve detects a hole, let us try to express the circulation along this curve, which we will call a \emph{global circulation}, as a sum of \emph{local circulations}, that is, circulations along small loops inside the region that it encloses. We will do this on two surfaces, a disk and a ring (Figure~\ref{fig:holes}). For the disk (top left panel), any two neighboring local circulations go along their shared boundary in opposite directions and cancel out there. Thus, if we fill the whole disk with local circulations, they cancel everywhere except on the boundary of the disk. In other words, the global circulation along the boundary is the net effect of the local circulations in the interior.\footnote{In essence, this argument is Stokes' theorem, $\oint_{\partial S}F\cdot d\ell=\int_S(\nabla\times F)\cdot dS$ for a vector field $F$ and a surface $S$ with boundary $\partial S$. The left-hand side is the global circulation along the curve $\partial S$, and the right-hand side is the sum of the local circulations, each captured by the curl of the field at the corresponding point.} For the ring (top middle panel), we can apply the same argument, obtaining not one but two global circulations, one along the outer boundary and another along the inner one. For the outer global circulation, the local circulations that generate it are inside the region that it encloses, as for the disk. For the inner one, however, we run into the problem that the local circulations generating it are \emph{outside} the region enclosed by it. We therefore cannot write this global circulation as a sum of local circulations in its interior. In fact, if we add local circulations on the surface of the ring, we can never obtain a global circulation along the inner boundary alone, it always comes together with one along the outer boundary circulating in the opposite direction. 

If we triangulate the two surfaces (bottom row of Figure~\ref{fig:holes}), the previous argument becomes the criterion behind Equation~\eqref{eq:H}, with the cycles around the triangles in the role of the local circulations. In the triangulated disk, the cycle along the boundary is the sum of the cycles around all the triangles, so it is in $\operatorname{im}B_2$ and is not harmonic. In contrast, in the triangulated ring, we can find a harmonic flow, namely the one that goes counterclockwise along both the inner and the outer cycle, with the same value on all their edges, and it has half this value on the edges between them. By construction, this flow cannot be expressed as a sum of triangle circulations, since the circulation of this flow around every triangle is zero. In fact, looking at the bottom row of Figure~\ref{fig:holes}, we can see a clear difference with the curl flow and the harmonic one: adding the circulations of all the triangles of the ring gives a flow that goes clockwise along the inner cycle (middle panel), while the harmonic flow goes counterclockwise along both cycles (right panel). Hence the holes that we found with $B_1$ and $B_2$ are analogous to the holes of the continuum.

So far, we have detected holes using $\operatorname{im}B_2\subseteq\ker B_1$ and finding zero-divergence flows that cannot be expressed as combinations of cycles around triangles. But we could use $\operatorname{im}B_1^T\subseteq\ker B_2^T$ and search for zero-curl flows that are not gradients. Since $\ker B_1=(\operatorname{im}B_1^T)^\perp$, the flows in $\ker B_2^T$ but not in $\im B_1^T$ are again the harmonic flows (the flows in $H$), so this second way of searching detects the same holes. Once again, this point of view has a counterpart in the continuum. For example, a classical way of showing that $\R^2\setminus\{0\}$ has a hole is to study the field $(-y,x)/(x^2+y^2)$. This field has zero curl, but its circulation along the unit circle is $2\pi$, so it cannot be the gradient of any function. On a simply connected domain, the Poincar\'e lemma states that every field with zero curl is a gradient. Hence the existence of this field implies that its domain has a hole, which in this case is the origin. The search for holes with fields instead of curves is called cohomology. Detecting holes with curves and detecting them with fields are in fact equivalent approaches, but proving it is hard (this equivalence is known as de Rham's theorem \cite{derham1931}). On a graph, however, the equivalence is an immediate consequence of the identity $\ker B_1\cap(\operatorname{im}B_2)^\perp=H=\ker B_2^T\cap(\operatorname{im}B_1^T)^\perp
$.

Let us now go back to Equation~\eqref{eq:H}. In Section~\ref{sec:subspaces}, we commented that $\R^m=\operatorname{im}B_1^T\oplus\ker B_1$. Moreover, we know that $\ker B_1=\operatorname{im}B_2\oplus H$. Hence
\begin{equation}\label{eq:hodge}
\R^m=\operatorname{im}B_1^T\oplus\operatorname{im}B_2\oplus H .
\end{equation}
In words, every flow is, in a unique way, the sum of a gradient, a flow of circulations around triangles and a harmonic flow. Equation~\eqref{eq:hodge} is the complete version of the Helmholtz decomposition of Section~\ref{sec:subspaces} \cite{eckmann1944,lim2020}. The harmonic part is nonzero only if the graph has holes, in the same way as the harmonic fields in the continuum.\footnote{The two incidence matrices also give the graph version of the relation between the Euler characteristic and the Betti numbers that we mentioned in Section~1. Let $\beta_0=\dim\ker B_1^T$ be the number of connected components, $\beta_1=\dim H$ the number of independent holes and $\beta_2=\dim\ker B_2$. If we apply the rank--nullity theorem to $B_1$ and $B_2$ and use that $\ker B_1=\operatorname{im}B_2\oplus H$, we obtain $n-m+k=\beta_0-\beta_1+\beta_2$. For example, for the house graph, $5-6+1=0=1-1+0$.}

To summarize, with two incidence matrices, we have understood the identity ``the boundary of a boundary is zero'', the holes of a space, the equivalence between homology and cohomology and the complete Helmholtz decomposition. In the curriculum, these results are spread over vector calculus in the first year and algebraic topology in the third year, if at all. They require chain complexes, singular homology and differential forms. On a graph, however, they are simple statements about the kernels and images of $B_1$ and $B_2$. We did not need group theory, quotients or analysis. Quite a payoff for a first course on linear algebra.

\section{Other branches of mathematics}\label{sec:other}

Analysis and topology are not the only subjects in which technical details obscure simple ideas. Many other courses of the curriculum start with a demanding definition and reach the central idea only after many technical steps. In this section, we will sketch, without complete proofs, three additional subjects in which graphs show the main ideas in a transparent way.

\subsection{Geometry}

On a curved surface, such as a sphere, we can still define a Laplace-like operator, called the Laplace--Beltrami operator. To define it, however, we need a Riemannian metric, which tells us how to measure lengths and angles at each point of the surface. The expression of the Laplace--Beltrami operator \cite{rosenberg1997} is notoriously scary,
\[
\Delta_g u=\frac{1}{\sqrt{\det g}}\sum_{i,j=1}^2\frac{\partial}{\partial x_i}\Bigl(\sqrt{\det g}\,(g^{-1})_{ij}\,\frac{\partial u}{\partial x_j}\Bigr),
\]
where $x_1,x_2$ are local coordinates on the surface and $g$ is the $2\times2$ matrix of the metric in these coordinates. Despite the intimidating expression of $\Delta_g$, the only new ingredient with respect to the usual Laplace operator is the metric. On a surface, we can always choose local coordinates in which $g$ is diagonal, with positive diagonal entries $g_1$ and $g_2$. In these coordinates, the formula becomes
\[
\Delta_g u=\frac{1}{\sqrt{g_1g_2}}\left[\frac{\partial}{\partial x_1}\Bigl(\sqrt{\tfrac{g_2}{g_1}}\,\frac{\partial u}{\partial x_1}\Bigr)+\frac{\partial}{\partial x_2}\Bigl(\sqrt{\tfrac{g_1}{g_2}}\,\frac{\partial u}{\partial x_2}\Bigr)\right].
\]
Up to the factor $1/\sqrt{g_1g_2}$, the operator takes the derivative of $u$ in each direction, multiplies it by a positive weight and then takes the divergence. This is the structure of the weighted Laplacian of a graph. Let us put a positive weight $w_e$ on each edge $e$, and let $W$ be the diagonal $m\times m$ matrix with the weights $w_e$ on its diagonal. The weighted Laplacian is then $BWB^T$, with the same gradient $B^T$ as in Section~\ref{sec:calculus}. Hence the weights give a gentle way to study how geometry changes the properties of the Laplacian, without inverse matrices, determinants or tensors.

\subsection{Probability and stochastic processes}

A central object of study in probability is the Markov process, a stochastic process in which the transition probabilities depend only on the current state. In particular, Brownian motion is used to model systems like the motion of a particle in a fluid or the fluctuations of stock prices. However, many important questions about Brownian motion can only be answered by solving partial differential equations. For example, the expected time that a Brownian particle needs to leave a region is the solution of a Poisson equation with Dirichlet boundary conditions.

When the number of states of a Markov process is finite, we can model it as a directed graph, with one vertex for each state and one edge for each possible transition, weighted by its probability \cite{levin2017}. Many properties of the process can then be studied through the properties of this graph. For instance, the process can go from any state to any other state if and only if the graph is strongly connected. Another example is that a strongly connected process is periodic (it can return to its starting state only at multiples of some integer larger than one) if and only if the lengths of all the cycles of the graph have a common divisor larger than one.

Brownian motion, in particular, can be modeled as a random walk on a connected \emph{undirected} graph in which each vertex $v$ moves to each of its neighbors with probability $1/d_v$. Interestingly, many expected times of this random walk can be computed, once again, with the Laplacian. For example, the expected time that the walk needs to leave a set of vertices is the solution of a linear system with the same matrix $L_{II}$ as the problem with zealots of Section~\ref{sec:diffusion}. Similarly, the expected time that the walk needs to go from a vertex $i$ to a vertex $j$ and back, called the commute time, is $\kappa(i,j)=2m\,(L^+_{ii}+L^+_{jj}-2L^+_{ij})$, where $L^+_{ij}$ are the entries of the pseudoinverse of $L$ \cite{chandra1996}. Hence quantities that in the continuum require the solution of partial differential equations are, on a graph, simple expressions in terms of the Laplacian.

\subsection{Optimization}

Linear programming is one of the most useful tools of applied mathematics, with applications in logistics,  economics, etc. A central result in the field is duality: every maximization problem has a dual minimization problem with the same optimal value. Students usually meet the dual problem as a formal recipe, which transposes the matrix of the problem and exchanges the roles of variables and constraints. The resulting problem seems to have no meaning of its own. Why should the maximum of one problem be equal to the minimum of a different one?

On a graph, both problems have a concrete meaning. Let us imagine that the edges are pipes, each with a capacity $c_e>0$, and that we want to send as much water as possible from a vertex $p$ to a vertex $q$. In the language of Section~\ref{sec:calculus}, we look for a flow $f$ with $|f_e|\le c_e$ on every edge and with zero divergence at every vertex other than $p$ and $q$. The dual problem looks for the cheapest way to block the water, that is, for a cut that separates $p$ from $q$ with the smallest total capacity. Every unit of water must cross every cut, so no flow can exceed the capacity of any cut. Less obviously, the largest flow is equal to the smallest cut. This result, known as the max-flow min-cut theorem  \cite{ford1956}, is a case of linear programming duality. Moreover, the maximization problem involves the divergence $B$, while the minimization problem involves the gradient $B^T$. Hence the transposition in the formal recipe of duality is simply a change of viewpoint from the divergence to the gradient.

\section{Discussion}\label{sec:discussion}

In this article, we have explored graphs as a setting in which to introduce key concepts of the undergraduate curriculum, such as the gradient, the divergence, the curl, the Laplacian, the diffusion equation, homology and the Helmholtz decomposition. Despite the diversity of these topics, most results were simple consequences of the properties of the kernels and images of incidence matrices. Hence, contrary to a common belief, understanding the central ideas of diffusion or homology, to name two examples, does not require an advanced knowledge of limits, function spaces or group theory. What students often perceive as difficult is not the central idea, but the analytic setting in which they first meet it. This does not mean that we should stop teaching the continuum. Instead, the continuum should come after the graph version, when the student already knows the central idea and can focus on the problems and subtleties of working over a continuum. In practice, a teacher needs only the linear algebra of a first course, so discrete spoilers can be introduced already in the second year of most curricula. A teacher can open each abstract topic with one or two lectures on its graph version, before introducing the continuous theory.

Several authors have already remarked that graphs can illustrate central results of their respective areas. For example, Doyle and Snell did so for random walks and potential theory \cite{doyle1984}, Lim for cohomology and Hodge theory \cite{lim2020}, Knill for the theorems of Green--Stokes, Gauss--Bonnet and Poincar\'e--Hopf \cite{knill2012}, and Grady and Polimeni for vector calculus and its applications to computational science \cite{grady2010}. Bamberg and Sternberg also introduced homology and cohomology to students of physics through electrical networks \cite{bamberg1990}. In contrast with these works, we propose graphs as a first setting for the curriculum as a whole. This proposal also challenges the widespread view that graphs should be taught primarily as models in applied mathematics. Graphs should also be leveraged in pure mathematics as a setting where a student can state a conjecture, test it on small examples and find the step where a proof fails.

A reader may object that we replace mathematics with a simplified imitation of it. Our answer is the maxim, usually attributed to Einstein, that everything should be made ``as simple as possible, but not simpler''. The graph versions are more intuitive than the continuous counterpart, but they are complete, with hypotheses, theorems and proofs. Only the setting is simpler; the mathematical content is the same.

Discrete spoilers are not a universal remedy and have their limits. First, our arguments are mathematical and historical. We have not measured the effect of the approach on students, which would require an empirical study in real courses. Second, the approach works only for theories whose central idea is evident in finite dimension. For example, measure theory has a central difficulty, the existence of sets without measure, that is not present in discrete settings. Other areas, like set theory, have no graph version at all, because their central ideas already require an infinite setting.

Spoilers are dreaded because knowing the ending can be fatal to a mediocre story. In great stories, however, the question is never what happens, but how. Everyone knows that Achilles will die, that Troy will fall and that Odysseus will need twenty years to return home, and yet the \emph{Iliad} and the \emph{Odyssey} still move readers almost three thousand years later. Everybody who watches \emph{Tristan und Isolde} knows that the lovers are doomed, but the four hours of musical longing make the final \emph{Liebestod} a cathartic experience. Mathematics is one of these stories. Telling a student, before any limit or function space, that a harmonic function is the average of its neighbors, or that a hole is a cycle that bounds nothing, does not spoil analysis or topology. It motivates the student to read on. We hope that our discrete spoilers leave students eager to learn not only which results are true, but why.

\section*{Acknowledgments}
The author thanks Karel Devriendt for insightful discussions. The author is also grateful to the podcast \textit{La vaca esf\'erica}, where he first heard the story of Euler's formula as told by Lakatos, and with it the view of mathematics as a sequence of conjectures, refutations and better proofs. The author acknowledges funding from PID2022-141802NB-I00 (BASIC), funded by MCIN/AEI/10.13039/501100011033 and by ``ERDF A way of making Europe'', and from Universidad Carlos III de Madrid through its Programa Propio de Investigaci\'on (Proyectos J\'ovenes PPIT2025, project IDISC, 2025/00800/001).

\end{document}